\documentclass[12pt, reqno]{amsart}
\usepackage{amsmath, amsthm, amscd, amsfonts, amssymb, graphicx, color}
\usepackage[bookmarksnumbered, colorlinks, plainpages]{hyperref}
\hypersetup{colorlinks=true,linkcolor=red, anchorcolor=green, citecolor=cyan, urlcolor=red, filecolor=magenta, pdftoolbar=true}

\newtheorem{theorem}{Theorem}[section]
\newtheorem{lemma}[theorem]{Lemma}
\newtheorem{proposition}[theorem]{Proposition}

\theoremstyle{definition}
\newtheorem{definition}[theorem]{Definition}

\theoremstyle{remark}

\numberwithin{equation}{section}

\begin{document}
\setcounter{page}{1}

\title[Eigenvalue decay of positive integral operators]{Eigenvalue decay of positive integral operators on compact two-point homogeneous spaces}

\author[M.H. Castro]{Mario Henrique de Castro}

\address{Faculdade de Matem\'atica, Universidade Federal de Uberl\^andia, Campus Santa M\^onica, Uberl\^andia, Brasil, CEP. 38.408-100}
\email{\textcolor[rgb]{0.00,0.00,0.84}{mariocastro@famat.ufu.br}}





\date{Received: xxxxxx; Revised: yyyyyy; Accepted: zzzzzz.
}

\begin{abstract}
%
%
We generalize and extend results on decay rates of singular values or eigenvalues of positive integral operators from unit spheres to two-point homogeneous spaces.\ The rates we present depend upon the order of the Laplace-Beltrami operator used to define the smoothness conditions on generating kernels, the Schatten class containing the integral operator generated by the derivative of the generating kernel and the dimension of the space.
\end{abstract} \maketitle

\section{Introduction}

Let $m\geq1$ be an integer and $\mathbb{M}$ be a compact two-point homogeneous space of dimension $m$.\ Such space is both a Riemannian $m$-manifold and a compact symmetric space of rank $1$.\ According to Wang \cite{wang}, they are the unit spheres $\mathbb{S}^m$, $m=1,2,\dots$; the real projective spaces $\mathbb{P}^m(\mathbb{R})$, $m=2,3,\dots$;  the complex projective spaces $\mathbb{P}^m(\mathbb{C})$, $m=4,6,\dots$; the quaternion projective spaces $\mathbb{P}^m(\mathbb{H})$, $m=8,12,16,\dots$; and Cayley's elliptic plane $\mathbb{P}^{16}$ of dimension $16$. More interesting information about these spaces can be found in \cite{cartan,gangolli,helgason1,helgason2}.

 In this paper, we will always consider $m\geq2$.\ Let $dx$ be the Riemannian measure on $\mathbb{M}$ and $L^2(\mathbb{M})$ the Hilbert space of all square-integrable complex functions on $\mathbb{M}$ endowed with the inner product
\begin{equation*}\langle
f,g\rangle_2:=\frac{1}{\sigma}\int_{\mathbb{M}}
f(x)\overline{g(x)}\,dx,\quad f,g \in L^2(\mathbb{M}),\end{equation*}
and the derived norm $||\cdot||_2$, the normalization constant being defined by $\sigma:=\int_{\mathbb{M}} \,dx$.

We will deal with integral operators defined by
\begin{equation}\label{eq1.1}\mathcal{K}(f)=\int_{\mathbb{M}} K(\cdot,
y)f(y)\,dy,\end{equation} in which the generating kernel $K\colon
\mathbb{M}\times \mathbb{M}\rightarrow\mathbb{C}$ is an element of
$L^2(\mathbb{M}\times \mathbb{M})$.\ In this case,
(\ref{eq1.1}) defines a compact operator on $L^2(\mathbb{M})$.
If $K$ is positive definite in the sense that
\begin{equation*}\int_{\mathbb{M}}\int_{\mathbb{M}} K(x,y)f(x)\overline{f(y)}\,dxdy\geq0,
\quad f \in L^2(\mathbb{M}),\end{equation*}
 then $\mathcal{K}$ becomes a
self-adjoint operator and the standard spectral theorem for compact and
self-adjoint operators is applicable and we can write
\begin{equation*}\mathcal{K}(f)=\sum_{n=0}^{\infty}\lambda_n(\mathcal{K})\langle f,
f_n\rangle_2 f_n, \quad f \in L^2(\mathbb{M}),\end{equation*} in which
$\{\lambda_n(\mathcal{K})\}$ is a sequence of nonnegative reals
(possibly finite) decreasing to 0 and $\{f_n\}$ is an $\langle
\cdot,\cdot\rangle_2$-orthonormal
basis of $L^2(\mathbb{M})$.\ The numbers $\lambda_n(\mathcal{K})$ are the
eigenvalues of $\mathcal{K}$ and the sequence
$\{\lambda_n(\mathcal{K})\}$ takes into account possible repetitions implied
by the algebraic multiplicity of each eigenvalue.\ The
positive definiteness of $K$ means nothing but the positivity of
the integral operator $\mathcal{K}$.\ Since it relates to the
inner product above, it is a common sense to call it $L^2$-{\em
positive definiteness}.

We observe that the addition of continuity to $K$ implies that
$\mathcal{K}$ is also {\em trace-class} (nuclear) (\cite{4,7,8}), that
is,
\begin{equation*}\sum_{f \in B}\langle \mathcal{K}^* \mathcal{K}( f),f\rangle_2^{1/2}
<\infty,\end{equation*} whenever $B$ is an orthonormal basis of
$L^2(\mathbb{M})$.\ In particular, it follows from Mercer's Theorem (\cite{mercer}) that
\begin{equation*}\sum_{n=1}^{\infty}\lambda_n(\mathcal{K})=\int_{\mathbb{M}}
K(x,x)\,dx<\infty,\end{equation*} and we can extract the most
elementary result on decay rates for the eigenvalues of such
operators, namely,
\begin{equation*}\lambda_n(\mathcal{K})=o(n^{-1}).\end{equation*}

If the integral operator $\mathcal{K}$ is compact but not
self-adjoint then decay rates for the singular values of the
operator becomes the focus.\ If $T$ is a compact operator on $L^2(\mathbb{M})$, its eigenvalues can be
ordered as $|\lambda_1(T)|\geq |\lambda_2(T)|\geq \cdots\geq 0$,
counting multiplicities (\cite{12}).\ The singular values of $T$ are, by
definition, the eigenvalues of the compact, positive and
self-adjoint operator $|T|:=(T^* T)^{1/2}$.\ The sequence
$\{s_n(T)\}$ of singular values of $T$ can also be ordered in a
decreasing manner, with repetitions being included according to
their multiplicities as eigenvalues of $|T|$.\ That being the case,
the classical Weyl's inequality (\cite[p.52]{8})
\begin{equation*}\Pi_{j=1}^n|\lambda_j(T)|\leq \Pi_{j=1}^n s_j(T), \quad n=1,2,\ldots,\end{equation*}
provides the convenient bridge between eigenvalues and singular
values.\ We remark that the inequality characterizing the
traceability of a compact non self-adjoint operator $T$ on
$L^2(\mathbb{M})$ reduces itself to
\begin{equation*}\sum_{n=1}^{\infty}s_{n}(T)< \infty\end{equation*}
and the elementary decay presented before becomes $s_n(\mathcal{K})=o(n^{-1})$.
Classical references on eigenvalues and singular values distribution
of compact operators on Banach spaces are \cite{12,19}.

The idea of nuclearity can be extended as follows.  For $p>0$ we say that a compact operator $T$ belongs to the {\em Schatten} $p$-{\em class} $\mathcal{S}_p$ if
\begin{equation*}
\sum_{n=1}^{\infty}{(s_n(T))^p}<\infty.
\end{equation*}
For $p\geq1$,  $\mathcal{S}_p$ is a Banach space endowed with the norm
$$\|T\|_p:=\left(\sum_{n=1}^{\infty}{(s_n(T))^p}\right)^{1/p}.$$
In particular, $\mathcal{S}_1$ and $\mathcal{S}_2$ coincide  respectively to the spaces of trace-class operators and of Hilbert-Schmidt operators.
Of course, if $T\in\mathcal{S}_p$ then its singular values satisfy $s_n(\mathcal{K})=o(n^{-1/p})$.

The problem of analyse the asymptotic behavior of
 $\{\lambda_n(\mathcal{K})\}$ or $\{s_n(\mathcal{K})\}$ under additional smoothness assumptions on
the kernel $K$ is the subject of this paper. Results of this very same nature can be found in \cite{mathcomp} (and references therein), where authors used Laplace-Beltrami differentiability as a condition of smoothness to the kernel $K$ and obtained some sharp results about how eigenvalues and singular values of $\mathcal{K}$ behave.\ The intention here is to invest in the very same issue by extending the setting from the sphere to two-point homogeneous spaces. However, in this more general setting we do not have a notion of devitative as the spherical Laplace-Beltrami derivative. Instead we use the Laplace-Beltrami operator to
define the smoothness condition we need.\ The connection between these two concepts on the sphere is relatively known and can be found in \cite{cmo}.  

The presentation of the paper is as follows.\ Section 2 contains
basic material about harmonic analysis in two-point homogeneous spaces and the
statement of the main result of the paper.\ In Section 3, we
state and prove some technical results related to Fourier expansions to be used in Section 4, where we present the proof for the main result
along with other pertinent information.

%

\section{Statement of the main results}

Two-point homogeneous spaces can be considered as the orbit of some compact subgroup $\mathfrak{H}$ of the orthogonal group $\mathfrak{G}$, i.e., $\mathbb{M}=\mathfrak{G}/\mathfrak{H}$.\ Let $e$ be the identity of $\mathfrak{G}$ and $\pi\colon\mathfrak{G}\to\mathfrak{G}/\mathfrak{H}$  the natural mapping. The {\em pole of     \,$\mathbb{M}$}, $o=\pi(e)$,  is invariant under all motions of $\mathfrak{H}$.\ Each one of these manifolds $\mathbb{M}$ has an invariant Riemannian metric $d(\cdot,\cdot)$ and a measure $dx$ induced by the normalized left Haar measure on $\mathfrak{G}$ which is invariant under the action of $\mathfrak{G}$. Also, these spaces admit essentially one invariant second order differential operator called {\em the Laplace-Beltrami operator} which we denote by $\Delta$. We suggest \cite{dai,kushpel2} and references therein for more information about this subject.

These spaces  have a similar geometry. For instance, all geodesics in one of them are closed and have the same length $2L$, in which
$L=\max\{d(x,y)\colon x,y\in\mathbb{M}\}$
is the diameter of $\mathfrak{G}/\mathfrak{H}$.\ A function on $\mathfrak{G}/\mathfrak{H}$ is invariant under the left action of $\mathfrak{H}$ on $\mathfrak{G}/\mathfrak{H}$ if and only if it depends only on the distance of its argument from the pole of $\mathbb{M}$.\ Let $\theta$ be the distance of a point from the pole.\ One can choose a geodesic polar coordinate system $(\theta,u)$, where $u$ is an angular parameter, in which the radial part of $\Delta$ can be written, up to a multiplicative constant, as
\begin{equation*}
\Delta_{\theta}=\frac{1}{(\sin{\lambda\theta})^{\sigma}(\sin{2\lambda\theta})^{\rho}}\frac{d}{d\theta}(\sin{\lambda\theta})^{\sigma}(\sin{2\lambda\theta})^{\rho}\frac{d}{d\theta},
\end{equation*}
in which
$$\begin{array}{|l|l|l|l|l|}
\hline
\mathbb{S}^m & \sigma=0 & \rho=m-1 & \lambda=\pi/2L & m=1,2,3,\dots; \\
\hline
\mathbb{P}^m(\mathbb{R}) & \sigma=0 & \rho=m-1 & \lambda=\pi/4L & m=2,3,4,\dots; \\
\hline
\mathbb{P}^m(\mathbb{C}) & \sigma=m-2 & \rho=1 & \lambda=\pi/2L & m=4,6,8,\dots; \\
\hline
\mathbb{P}^m(\mathbb{H}) & \sigma=m-4 & \rho=3 & \lambda=\pi/2L & m=8,12,\dots;\\
\hline
\mathbb{P}^{16}(Cay) & \sigma=8 & \rho=7 & \lambda=\pi/2L & m=16.\\
\hline
\end{array}$$
Furthermore, the change of variables $x=\cos{2\lambda\theta}$ gives us
\begin{equation*}
\Delta_x=(1-x)^{-\alpha}(1+x)^{-\beta}\frac{d}{dx}(1-x)^{1+\alpha}(1+x)^{1+\beta}\frac{d}{dx},
\end{equation*}
with $\alpha=(\sigma+\rho-1)/2=(m-2)/2$ and $\beta=(\rho-1)/2$.

We will write $\mathcal{B}=-\Delta_x$ and also call it Laplace-Beltrami operator on $\mathbb{M}$.\ Let denote $\mathcal{B}^r$ the $r$-th power of $\mathcal{B}$, $r=0,1,2,\dots$.\ The {\em Sobolev space of order $r$} constructed from $\mathcal{B}$ is defined as in \cite[p.37]{lions} and \cite{platonov1} by
\begin{equation*}
W_2^r(\mathbb{M}):=\{f\in L^2(\mathbb{M})\colon \mathcal{B}^jf\in L^2(\mathbb{M}),\, j=1,2,\dots,r\},
\end{equation*}
and the Laplace-Beltrami operator and its powers satisfy
\begin{equation}\label{eq3.2}\langle \mathcal{B}^r f, g \rangle_2 = \langle f, \mathcal{B}^r g\rangle_2,\quad  f,g \in W_2^r(\mathbb{M}).\end{equation}
The action of the Laplace-Beltrami operator on kernels is done
separately: we keep one variable fixed and differentiate with
respect to the other.\ We denote $K_{0,r}(x, y)$ to indicate
the $r$-th order of $\mathcal{B}$ acting on the kernel $K$ with respect to the second
variable $y$ (we do not differentiate with respect to the first
variable $x$). The integral operator associated
with $K_{0,r}$ will be written as $\mathcal{K}_{0,r}$.

\begin{definition} A kernel $K\in L^2(\mathbb{M}\times\mathbb{M})$ belongs to $W_2^r$ if $K(x,\cdot)\in W_2^r$, $x\in \mathbb{M}$ a.e..\end{definition}

\noindent We are able to state the main result we intend to prove here.

\begin{theorem}\label{teo2.3}Let $K\in L^2(\mathbb{M}\times
\mathbb{M})$ be a $L^2$-positive
definite kernel  satisfying $K(x,\cdot)\in W_2^r(\mathbb{M})$, $x\in\mathbb{M}$ a.e..\ If $\mathcal{K}_{0,r}\in \mathcal{S}_p$
then
\begin{equation*}\lambda_n(\mathcal{K})=o(n^{(-1/p)-(2r/m)}).\end{equation*}\end{theorem}

\noindent{\bf Remark 2.3.} This result unifies and generalizes Theorems 2.4 and 2.5 in \cite{mathcomp}. Moreover, Theorem 2.3 of \cite{mathcomp} can also be generalized to homogeneous spaces but since the proofs are very similar then we do not state e prove it here.

\section{Fourier analysis for functions on $\mathbb{M}$}

The Hilbert space $L^2(\mathbb{M})$ can be decomposed as $L^2(\mathbb{M})=\oplus_{n=0}^{\infty}\mathcal{H}_n^m$, where $\mathcal{H}_n^m$ is the eigenspace of $\mathcal{B}$ with respect to the eigenvalue $\lambda_n(\mathcal{B}):=n(n+\alpha+\beta+1)$ ($\lambda_0(\mathcal{B})=1$).\ The elements of  $\mathcal{H}_n^m$ are the well-known Jacobi polynomials $P_n^{(\alpha,\beta)}$.
Each $\mathcal{H}_n^m$ has a finite dimension given by the formula
\begin{equation*}\label{dn}
d_n^m=d_n^m(\mathbb{M})=\frac{\Gamma(\beta+1)(2n+\alpha+\beta+1)\Gamma(n+\alpha+1)\Gamma(n+\alpha+\beta+1)}{\Gamma(\alpha+1)\Gamma(\alpha+\beta+2)\Gamma(n+1)\Gamma(n+\beta+1)},
\end{equation*}
for all $n\in\mathbb{N}$, if $\mathbb{M}\in\{\mathbb{S}^m,\mathbb{P}^m(\mathbb{C}),\mathbb{P}^m(\mathbb{H}),\mathbb{P}^{16}(Cay)\}$ and for  $n$ even if $\mathbb{M}=\mathbb{P}^m(\mathbb{R})$.  Otherwise, for n odd, $d_n^m(\mathbb{P}^m(\mathbb{R}))=0$.  
Let denote $\mathcal{T}_n^m=\oplus_{k=0}^n\mathcal{H}_n^m$ and $\tau_n^m=\dim{\mathcal{T}_n^m}$.
Precisely, for $\mathbb{S}^m,\mathbb{P}^m(\mathbb{C}),\mathbb{P}^m(\mathbb{H})$, and $\mathbb{P}^{16}(Cay)$ we have the explicit expression
\begin{equation}\label{tau}
\tau_n^m=\frac{\Gamma(\beta+1)\Gamma(n+\alpha+\beta+2)\Gamma(n+\alpha+2)}{\Gamma(\alpha+\beta+2)\Gamma(\alpha+2)\Gamma(n+\beta+1)\Gamma(n+1)},
\end{equation}for all $n\in\mathbb{N}$.
It follows that
\begin{equation}\label{dlimit}d_n^m=O(n^{m-1}), \quad\mbox{as } n\to\infty.
\end{equation} 	
\begin{equation}\label{taulimit}\tau_n^m=O(n^m), \quad\mbox{as } n\to\infty.
\end{equation}

Let $\{Y_{n,k}: k=1,\dots,d_n^m\}$ be an orthonormal basis of $\mathcal{H}_n^m$.\ Each $f\in L^2(\mathbb{M})$ has a Fourier expansion
$$f=\sum_{n=0}^{\infty}\sum_{k=1}^{d_n^m}\, c_{n,k}(f)\, Y_{n,k},$$
in which $c_{n,k}(f)=\langle f,Y_{n,k}\rangle_2$.

If $r > 0$, a function $g \in L^2(S^m)$ is called the {\it fractional derivative of order $r$ of $f$} whenever the Fourier series of $g$ has the form
$$g=\sum_{n=1}^{\infty}\sum_{k=1}^{d_n^m}n^r(n+\alpha+\beta-1)^r\,c_{n,k}(f)\,Y_{n,k}.$$
Thus, the fractional derivative of order $r$ of a function $f\in W_2^r(\mathbb{M})$ corresponds to $\mathcal{B}^rf$ whenever $r$ is a positive integer.

The $r$-th fractional integral ($r>0$) of $f\in L^2(\mathbb{M})$ is the $L^2(\mathbb{M})$ element
\begin{equation}\label{int op}J^rf=c_{0,1}(f)+\sum_{n=1}^{\infty}\sum_{k=1}^{d_n^m}n^{-r}(n+\alpha+\beta-1)^{-r}\,c_{n,k}(f)\,Y_{n,k}.\end{equation}

\begin{proposition} Let $r>0$.\ The {\it $r$-th fractional integral operator} $J^r:L^2(\mathbb{M})\to L^2(\mathbb{M})$ defined by \eqref{int op} is a compact operator.\end{proposition}

\begin{proof} It is enough to observe that $J^r$ is linear and can be approximated  by some sequence of linear finite rank operators in the space of bounded operators on $L^2(\mathbb{M})$.\end{proof}

\begin{proposition} If $r$ be a positive integer then $J^r(L^2(\mathbb{M}))\subset W_2^r(\mathbb{M})$.
 \end{proposition}

\begin{proposition}\label{integral_inverse} If $r$ is a positive integer 
and $f \in \oplus_{n=1}^{\infty}\mathcal{H}_n^m$ then $\mathcal{B}^r J^r f=f$.\end{proposition}

The singular values of the $J^r$ are given by $s_0(J^r)=1$ and
\begin{equation*}s_n(J^r)=\lambda_n(\mathcal{B}^r)^{-1}=n^{-r}(n+\alpha+\beta+1)^{-r},\quad n=1,2,\dots.
\end{equation*}
  They are ordered in
accordance with the spectral theorem for compact
operators.\ In other words, we assume they are listed in decreasing order counting the repetitions.
As so, if $\mathbb{M}\in\{\mathbb{S}^m,\mathbb{P}^m(\mathbb{C}),\mathbb{P}^m(\mathbb{H}),\mathbb{P}^{16}(Cay)\}$ we may think the sequence
                                     $\{s_n(J^r)\}$ is
block ordered in such a way that the first block contains the
singular value $s_0(J^r)=1$ and the $(n+1)$-th block ($n\geq 1$) contains $d_n^m$
entries equal to $n^{-r}(n+\alpha+\beta+1)^{-r}$.\ For future reference,
we notice that the first entry in the $(n+1)$-th block corresponds
to the index
\begin{equation*}d_0^m+d_1^m+\cdots +d_{n-1}^m+1=\tau_{n-1}^m+1.\end{equation*}
As for the last one, it corresponds to
\begin{equation}\label{tauene}d_0^m+d_1^m+\cdots +d_{n-1}^m+d_n^m=\tau_{n}^m.\end{equation}

In the next lemmas we detach technical inequalities to be
used ahead.\ The first one improves  \eqref{taulimit}.

\begin{lemma}\label{lem4.4} If $\mathbb{M}\in\{\mathbb{S}^m,\mathbb{P}^m(\mathbb{C}),\mathbb{P}^m(\mathbb{H}),\mathbb{P}^{16}(Cay)\}$ then there exists an integer $\delta(m)\geq
1$ such that
\begin{equation*}\tau_n^m\leq 2n^m,\quad n\geq\delta(m).\end{equation*}
\end{lemma}

\begin{proof} We keep \eqref{taulimit} in mind and develop each case using \eqref{tau} in order to obtain  a polynomial expression to $\tau_n^m$ .\ For $\mathbb{M}=\mathbb{S}^m$ we know $\alpha=\beta=(m-2)/2$.\ Consequently there is $\delta(\mathbb{S}^m)>0$ such that
\begin{equation*}\label{tauesfera}\tau_n^m=\frac{2n^m}{m!}\left(1+\frac{c_{1}^{(1)}}{n}+\frac{c_2^{(1)}}{n^2}+\cdots+\frac{c_m^{(1)}}{n^m}\right)\leq2n^m, \quad n\geq\delta(\mathbb{S}^m),
\end{equation*} where $c_1^{(1)},\dots,c_m^{(1)}$ do not depend upon $n$.

For $\mathbb{M}=\mathbb{P}^m(\mathbb{C})$ we know $\alpha=(m-2)/2$ and $\beta=0$.\ As so, there is $\delta(\mathbb{P}^m(\mathbb{C}))>0$ such that
\begin{equation*}\label{tauprojcomp} \tau_n^m=\left[\frac{n^{m/2}}{(m/2)!}\left(1+\frac{c_1^{(2)}}{n}+\cdots+\frac{c_m^{(2)}}{n^{m/2}}\right)\right]^2\leq2n^m,\quad n\geq\delta(\mathbb{P}^m(\mathbb{C})),
\end{equation*} in which $c_1^{(2)},\dots,c_m^{(2)}$ do not depend upon $n$.

If $\mathbb{M}=\mathbb{P}^m(\mathbb{H})$ then $\alpha=(m-2)/2$ and $\beta=1$.\ Thus, there is $\delta(\mathbb{P}^m(\mathbb{H}))>0$ such that
\begin{equation*}\label{tauprojquat} \tau_n^m=\frac{(n+1+m/2)}{(n+1)(1+m/2)}\left(\frac{(n+m/2)!}{n!(m/2)!}\right)^2\leq2n^m,\quad n\geq\delta(\mathbb{P}^m(\mathbb{H})).
\end{equation*}

If $\mathbb{M}=\mathbb{P}^{16}(Cay)$ then $\alpha=(m-2)/2$, $\beta=3$, and there is $\delta(\mathbb{P}^{16}(Cay))>0$ such that
\begin{equation*}\label{taucay} \tau_n^{16}=\frac{(n+12)(n+11)(n+10)(n+9)}{1980(n+4)(n+3)(n+2)(n+1)}\left(\frac{(n+8)!}{(n)!8!}\right)^2\leq2n^{16},
\end{equation*} since $n\geq\delta(\mathbb{P}^{16}(Cay))$.

To conclude we  define $\delta(m)=\max\{\delta(\mathbb{S}^m),\delta(\mathbb{P}^m(\mathbb{C})),\delta(\mathbb{P}^m(\mathbb{H})),\delta(\mathbb{P}^{16}(Cay))\}$.
\end{proof}

\begin{lemma}\label{lem4.5}If $m$ is an integer at least 2 then
\begin{equation*}(n+1)^{m}-(n^{m}+1)+1\leq m2^{m-1}n^{m-1},\quad n\geq 1.\end{equation*}\end{lemma}

\begin{proof} It suffices to apply the mean value theorem
to the function $x^m$ on the interval $[n,n+1]$ and estimate the
resulting formula conveniently.\end{proof}

\section{Proof of the main result}

This section contains a proof for Theorem \ref{teo2.3}.\ It depends upon
some general properties of compact operators and their singular
values which we now describe in a form adapted to our needs.\ They
can be found in standard references on operator theory such as
\cite{7,8,12,19} and depend on the ordering of eigenvalues and singular
values as previously mentioned.

\begin{lemma}\label{lem4.1}Let $T$ be a compact operator on
$L^2(\mathbb{M})$.\ The following assertions hold:\\
$(i)$ If $T$ is self-adjoint then
\begin{equation*}s_{n}(T)=|\lambda_{n}(T)|,\quad n=1,2,\ldots;\end{equation*}
$(ii)$ If $A$ is a bounded operator on $L^2(\mathbb{M})$ then both, $AT$
and $TA$, are compact.\ In addition,
\begin{equation*}\max\{s_{n}(AT), s_{n}(TA)\} \leq \|A\|\, s_{n}(T), \quad
n=1,2,\ldots;\end{equation*}
$(iii)$ If $A$ is a linear operator on
$L^2(\mathbb{M})$ of rank at most $l$, then
\begin{equation*}s_{n+l}(T) \leq s_{n}(T+A), \quad n=1,2,\ldots;\end{equation*}
 $(iv)$ If $A$ is a compact
operator on $L^2(\mathbb{M})$ then
\begin{equation*}s_{n+k-1}(AT) \leq
s_{n}(A)s_{k}(T),\quad n,k=1,2,\ldots.\end{equation*}\end{lemma}

The following additional lemma regarding the singular values of an
integral operator generated by a square-integrable kernel is proved
in \cite[p.40]{12}.

\begin{lemma}\label{lem4.2}If $K\in L^2(\mathbb{M}\times \mathbb{M})$ then
\begin{equation*}\sum_{n=1}^{\infty}s_{n}^2(\mathcal{K})=\|K\|^2_{2}.\end{equation*}\end{lemma}

The key idea behind the proof of the main result previously stated
resides in the following estimation for the singular values of
$\mathcal{K}$, which holds when $K$ is smooth enough.

\begin{lemma}\label{lem4.3}Let $K$ be an element of $W_2^r(\mathbb{M})$.\ If
$\mathcal{K}_{0,r}$ is bounded then
\begin{equation*}s_{n+1}(\mathcal{K})\leq s_n(\mathcal{K}_{0,r}J^r),\quad n=1,2,\dots.\end{equation*}\end{lemma}

\begin{proof} Consider the orthogonal projection $Q$ of
$L^2(\mathbb{M})$ onto $\oplus_{\ell=1}^{\infty}\mathcal{H}_{\ell}^{m+1}$.\
Since $I-Q$ is a projection onto the orthogonal complement of
$\oplus_{\ell=1}^{\infty}\mathcal{H}_{\ell}^{m+1}$ then
$\mathcal{K}-\mathcal{K}Q$ is an operator on $L^2(\mathbb{M})$ of rank at
most 1.\ Using Lemma \ref{lem4.1}-$(iii)$, we may deduce that
\begin{equation}\label{eq4.7}s_{n+1}(\mathcal{K})\leq
s_{n}(\mathcal{K}-\mathcal{K}(I-Q))=s_{n}(\mathcal{K}Q), \quad
n=1,2,\ldots.\end{equation} To proceed, we need a convenient
decomposition for $\mathcal{K}Q$.\ Looking at the action of
$\mathcal{K}Q$ on a generic element $f$ from $L^2(\mathbb{M})$ and using
Proposition \ref{integral_inverse} we see that
\begin{equation*}\mathcal{K}Q(f)=\int_{\mathbb{M}}K(\cdot,y)Qf(y)\,d\sigma_m(y)=\int_{\mathbb{M}}K(\cdot,y)
\mathcal{B}^rJ^rQf(y)\,d\sigma_m(y).\end{equation*} Since $K\in W_2^r(\mathbb{M})$,
we employ (\ref{eq3.2}) to obtain
\begin{equation*}\mathcal{K}Q(f)=\int_{\mathbb{M}} K_{0,r}(\cdot,y)J^r(Qf)(y)\,d\sigma_m(y)=\mathcal{K}_{0,r}J^rQ(f),\end{equation*}
that is, $\mathcal{K}Q=\mathcal{K}_{0,r}J^r Q$.\ Now, assuming
$\mathcal{K}_{0,r}$ is bounded, we can apply (\ref{eq4.7}) and Lemma
\ref{lem4.1}-$(ii)$ to see that
\begin{equation*}
s_{n+1}(\mathcal{K})\leq s_{n}(\mathcal{K}Q)\leq
\|Q\|s_{n}(\mathcal{K}_{0,r}J^r)\leq
s_{n}(\mathcal{K}_{0,r}J^r),\quad n=1,2,\ldots.\end{equation*} The
proof is complete.\end{proof}

The following technical result is borrowed from \cite{13}. An elementary proof of such result can be found in \cite{azevedo}.

\begin{lemma}\label{lem4.7}Let $\{a_n\}$ be a decreasing
sequence of positive real numbers.\ If the series
$\sum_{n=1}^{\infty}n^{\alpha}a_n^{\beta}$ is convergent for some
positive constants $\alpha$ and $\beta$ then
$a_n=o(n^{-(\alpha+1)/\beta})$.\end{lemma}

We now proceed to the proof of the main result in the
paper.\\

\begin{proof}[{\bf Proof of Theorem \ref{teo2.3}}] We perform the demonstration in three steps. First, we assume $\mathcal{K}_{0,r}$
belongs to $\mathcal{S}_p$ and show
$$\sum_{n=1}^{\infty}n^{2rp+m-1}(\lambda_{n^m}(\mathcal{K}))^p <\infty.$$
Second, we prove
$$\sum_{n=1}^{\infty}n^{2rp/m}(\lambda_n(\mathcal{K}))^p<\infty.$$
Finally, we apply Lemma \ref{lem4.7} to conclude that
$$\lim_{n\to\infty}n^{{1\over p}+{2r\over m}}\lambda_n(\mathcal{K})=0.$$
Combining Lemma \ref{lem4.3} with Lemma \ref{lem4.1}-$(iv)$ we can deduce the
inequalities
\begin{equation*}s_{n+k}(\mathcal{K})\leq s_{n+k-1}(\mathcal{K}_{0,r}J^r)\leq s_{k}(\mathcal{K}_{0,r})s_{n}(J^r),\quad n,k=1,2,\dots,
\end{equation*}
Since the sequence of eigenvalues of $\mathcal{K}$ is non increasing, it follows from Lemma \ref{lem4.1}-(i) that
\begin{equation*}
\lambda_{\tau_n+k}(\mathcal{K})\leq s_{k}(\mathcal{K}_{0,r})s_{\tau_n}(J^r)=
s_k(\mathcal{K}_{0,r})\ n^{-r}(n+\alpha+\beta+1)^{-r},\quad n,k=1,2,\dots,
\end{equation*} where the last equality is a consequence of \eqref{tauene}.
Thus
\begin{equation*}n^{2r}\lambda_{\tau_n+k}(\mathcal{K})\leq n^r (n+\alpha+\beta+1)^r\lambda_{\tau_n+k}(\mathcal{K})\leq s_{k}(\mathcal{K}_{0,r}),\ n,k=1,2,\ldots,
\end{equation*}
which implies
\begin{equation*}
n^{2r p}(\lambda_{\tau_n+k}(\mathcal{K}))^{p}\leq (s_{k}(\mathcal{K}_{0,r}))^{ p},\quad n,k=1,2,\dots.
\end{equation*}
Since $\mathcal{K}_{0,r}\in \mathcal{S}_p$, by adding on $k$ and $n$ leads to
\begin{eqnarray*}
\sum_{n=1}^{\infty}n^{2rp}\sum_{k=\tau_{n-1}+1}^{\tau_n}(\lambda_{\tau_n+k}(\mathcal{K}))^p
&\leq& \sum_{n=1}^{\infty}\sum_{k=\tau_{n-1}+1}^{\tau_n}(s_{k}(\mathcal{K}_{0,r}))^p\\
&\leq& \sum_{n=1}^{\infty}(s_{n}(\mathcal{K}_{0,r}))^p <\infty.
\end{eqnarray*}
To proceed, we apply
Lemma \ref{lem4.4} to select a constant $\delta=\delta(m)\geq 1$ such that
\begin{equation*}2\tau_n^m\leq 2^2n^{m} \leq (2n)^{m},\quad
n\geq \delta.\end{equation*}
As long as \eqref{dlimit} gives us a constant $c>0$ such that
\begin{equation*}
n^{m-1}\leq c\ d_n^m, \quad n\geq \beta=\beta(m),
\end{equation*}
and
$\{\lambda_n(\mathcal{K})\}$ does not increases, choosing $\gamma=\max\{\delta(m),\beta(m)\}$, we now see that
\begin{align*}
\sum_{n\geq\gamma}(2n)^{2rp+m-1}(\lambda_{(2n)^{m}}(\mathcal{K}))^p
&\leq c\sum_{n\geq\gamma}(2n)^{2rp}d_n^m(\lambda_{(2n)^{m}}(\mathcal{K}))^p\\
&\leq  2^{2rp}c\sum_{n\geq\gamma}n^{2rp}\sum_{k=\tau_{n-1}+1}^{\tau_n}(\lambda_{(2n)^{m}}(\mathcal{K}))^p\\
&\leq c_1\sum_{n\geq\gamma}n^{2rp}\sum_{k=\tau_{n-1}+1}^{\tau_n}(\lambda_{2\tau_n}(\mathcal{K}))^p\\
&\leq c_1\sum_{n\geq\gamma}n^{2rp}\sum_{k=\tau_{n-1}+1}^{\tau_n}(\lambda_{\tau_n+k}(\mathcal{K}))^p
<\infty.\end{align*}
Moreover, we show in the same way there is $c_2=c_2(r,p,m)>0$ such that
\begin{align*} \sum_{n\geq\gamma}(2n+1)^{2rp+m-1}(\lambda_{(2n+1)^{m}}(\mathcal{K}))^p
&\leq c_2\sum_{n\geq\gamma}n^{2rp}\sum_{k=\tau_{n-1}+1}^{\tau_n}(\lambda_{\tau_n+k}(\mathcal{K}))^p
<\infty.\end{align*}
Hence, we conclude the first step showing that
\begin{equation*}
\sum_{n=1}^{\infty}n^{2rp+m-1}(\lambda_{n^m}(\mathcal{K}))^{p}<\infty.
\end{equation*}
The second step starts with the notice that
\begin{eqnarray*}\sum_{n=1}^{\infty} n^{2rp\over m}(\lambda_n(\mathcal{K}))^p
 &=&\sum_{n=1}^{\infty} \sum_{k=0}^{(n+1)^{m}\mbox{-}(n^{m}+1)}(n^{m}+k)^{2rp\over m}
(\lambda_{n^{m}+k}(\mathcal{K}))^p\\
&\leq &\sum_{n=1}^{\infty}\sum_{k=0}^{(n+1)^{m}\mbox{-}(n^{m}+1)}(n^{m}+(n+1)^m-n^m-1)^{2rp\over m}
(\lambda_{n^{m}+k}(\mathcal{K}))^p.
\end{eqnarray*}
With an application of Lemma \ref{lem4.5} we obtain
\begin{eqnarray*}\sum_{n=1}^{\infty} n^{2rp\over m}(\lambda_n(\mathcal{K}))^p
&\leq &\sum_{n=1}^{\infty}\sum_{k=0}^{(n+1)^{m}\mbox{-}(n^{m}+1)}\left[(2n)^m\right]^{2rp\over m}
(\lambda_{n^{m}+k}(\mathcal{K}))^p\\
&\leq & 2^{2rp}\sum_{n=1}^{\infty}n^{2rp}(\lambda_{n^{m}}(\mathcal{K}))^p\sum_{k=0}^{(n+1)^{m}\mbox{-}(n^{m}+1)}1\\
&= & 2^{2rp}\sum_{n=1}^{\infty}n^{2rp}(\lambda_{n^{m}}(\mathcal{K}))^p[(n+1)^{m}-n^{m}],
\end{eqnarray*}
from which we can find a constant $c_3=c_3(r,p,m)>0$ so that
\begin{eqnarray*}\sum_{n=1}^{\infty} n^{2rp\over m}(\lambda_n(\mathcal{K}))^p&\leq & c_3\sum_{n=1}^{\infty}n^{2rp}(\lambda_{n^{m}}(\mathcal{K}))^p\,n^{m-1}\\
&\leq& c_3\sum_{n=1}^{\infty}n^{2rp+m-1}(\lambda_{n^m}(\mathcal{K}))^{p}<\infty.
 \end{eqnarray*}
Finally, Lemma \ref{lem4.7} is applied to give
$$\lim_{n\to\infty}n^{{1\over p}+{2r\over m}}\lambda_n(\mathcal{K})=0$$
and complete the proof.
\end{proof}

\noindent{\bf Remark 4.5.} The proof we  performed here does not include the case $\mathbb{M}=\mathbb{P}^m(\mathbb{R})$. However, Theorem \ref{teo2.3} works also in this case. Indeed, functions on  $\mathbb{P}^m(\mathbb{R})$ can be seen as even functions on $\mathbb{S}^m$, so $L^2(\mathbb{P}^m(\mathbb{R}))$ can be identified to $\oplus_n\mathcal{H}_{2n}^m(\mathbb{S}^m)$ and the case $\mathbb{M}=\mathbb{P}^m(\mathbb{R})$ follows directly from the case $\mathbb{M}=\mathbb{S}^m$.\\

\noindent{\bf Remark 4.6.} The decay rate obtained was proved to be optimal in \cite{mathcomp} for the case $\mathbb{M}=\mathbb{S}^m$, with $p=1$.\\

\noindent{\bf Acknowledgement.} I thank Professor Valdir A. Menegatto for suggesting the theme.\\
This work was partially supported by FAPEMIG, Grant APQ-00474-14; and CNPq, Grant 475320/2013-1.

\bibliographystyle{amsplain}

\end{document}